\documentclass[reqno]{amsart}
\usepackage{latexsym}
\usepackage{amscd}
\usepackage[dvips]{graphicx}
\newtheorem{thm}{Theorem}

\newenvironment{pf}{\begin{proof}}{\end{proof}}

\newtheorem{dfn}{Definition}[section]

\newtheorem{prp}{Proposition}[section]

\newtheorem{lma}{Lemma}[section]

\newtheorem{rem}{Remark}[section]
\newenvironment{rmk}{\begin{rem}\rm}{\end{rem}}

\newcommand{\tR}{{\mathbb{R}}}

\newcommand{\tZ}{{\mathbb{Z}}}
\newcommand{\tF}{{\mathcal{F}}}
\newcommand{\Si}{\Sigma}

\newcommand{\wt}{\widetilde}

\newcommand{\lk}{\operatorname{lk}}

\newcommand{\Imm}{{\mathbf {Imm}}}
\newcommand{\Emb}{{\mathbf  {Emb}}}
\newcommand{\krn}{\operatorname{ker}}
\newcommand{\cokrn}{\operatorname{coker}}

\newcommand{\St}{\operatorname{St}}

\newcommand{\inr}{\operatorname{int}}

\newcommand{\Or}{\operatorname{Or}}
\newcommand{\Tor}{\operatorname{Tor}}
\renewcommand{\mod}{\operatorname{mod}}
\newcommand{\skarp}{\operatorname{\ensuremath{\sharp}}}

\begin{document}
\title
[Invariants of generic immersions and Bernoulli numbers]
{Invariants of generic immersions\\ and Bernoulli numbers}
\author{Tobias Ekholm}
\address{Department of Mathematics, Uppsala University, S-751 06 Uppsala,
Sweden}
\email{tobias@math.uu.se}
\date{}
\begin{abstract}
First order invariants of generic immersions of manifolds of dimension
$nm-1$ into manifolds of dimension $n(m+1)-1$, $m,n>1$ are
constructed using the geometry of self-intersections. 
The range of one of these invariants is related to Bernoulli numbers. 

As by-products some geometrically defined invariants of regular
homotopy are found.  
\end{abstract}
\maketitle

\section{Introduction}
An immersion of a smooth manifold $M$ into a smooth manifold $W$ is
a smooth map with everywhere injective differential. Two immersions
are regularly homotopic if they can be connected by a continuous
1-parameter family of immersions. 

An immersion is generic if all its self-intersections are
transversal. In the space $\tF$ of immersions $M\to W$, generic 
immersions form an open dense subspace. Its complement is 
{\em the discriminant hypersurface}, $\Si\subset\tF$.
Two generic immersions belong to the same path component of $\tF-\Si$
if they can be connected by a regular homotopy, which at
each instance is a generic immersion. We shall consider
the classification of generic immersions up to regular homotopy
through generic immersions. It is similar to the
classification of embeddings up to diffeotopy (knot theory). 
In both cases, all topological properties of equivalent maps are the
same. 

An {\em invariant of generic immersions} is a function on $\tF-\Si$
which is locally constant. The value of such a function along a 
path in $\tF$ jumps at intersections with $\Si$. Invariants may be
classified according to the complexity of their jumps. The most basic
invariants in this classification are called 
{\em first order invariants} (see Section ~\ref{mainres}).  

In \cite{A}, Arnold studies generic regular plane curves (i.e. generic
immersions $S^1\to\tR^2$). He finds three 
first order invariants $J^+$, $J^-$, and
$\St$. In \cite{E2}, the author considers the case $S^3\to\tR^5$. Two
first order 
invariants $J$ and $L$ are found. In both these cases, the only
self-intersections of generic immersions are transversal double points 
and in generic 1-parameter families there appear isolated
instances of self-tangencies and triple points. The values of the
invariants $J^{\pm}$ and $J$ change at instances of self-tangency and 
remain constant at instances of triple points. The invariants $\St$
and $L$ change at triple points and remain constant at
self-tangencies. 

In this paper we shall consider high-dimensional analogs of these
invariants. The most straightforward generalizations arise for 
immersions $M^{2n-1}\to W^{3n-1}$, where self-tangencies and
triple points are the only degeneracies in generic 1-parameter families.
The above range of dimensions is included in a
2-parameter family, $M^{nm-1}\to W^{n(m+1)-1}$, $m,n>1$, where
generic immersions do not have $k$-fold self-intersection points if
$k>m$ and in generic 1-parameter families there appear isolated
instances of $(m+1)$-fold self-intersection. Under these circumstances, 
we find first order invariants of generic immersions. 
(For precise statements, see Theorem
~\ref{thmL} and Theorem ~\ref{thmJ} in Section ~\ref{mainres}, where
the main results of this paper are formulated.) 

In particular, if $n$ is
even and $M^{nm-1}$ is orientable and satisfies a certain homology
condition (see Theorem ~\ref{thmL} (b)) then there exists an
integer-valued invariant $L$ of generic immersions 
$M^{nm-1}\to\tR^{n(m+1)-1}$ which is an analog of Arnold's $\St$: It
changes under instances of $(m+1)$-fold self-intersection, does not 
change under other degeneracies which appear in generic 1-parameter
families, and is additive under connected sum. The value of $L$
at a generic immersion $f$ is the linking number of a
copy of the set of $m$-fold self-intersection points of $f$, shifted in 
a special way, with $f(M)$ in $\tR^{n(m+1)-1}$. (See Definition
~\ref{dfnL}.)     

For immersions of odd-dimensional spheres in codimension two there
appear restrictions on the possible values of $L$ (see Theorem
~\ref{thmdivL}). This phenomenon is especially interesting in the case
of immersions $S^{4j-1}\to\tR^{4j+1}$. Using a result of Hughes and
Melvin \cite{HM} we prove that the range of $L$ is 
related to Bernoulli numbers: For any 
generic immersion $f\colon S^{4j-1}\to\tR^{4j+1}$, $L(f)$ is divisible 
by $p^r$, where $p$ is a prime and $r$ is an integer such that
$p^{r+k}$ divides $2j+1$ and $p^{k+1}$ does not divide $\mu_j$ for
some integer $k$. Here $\mu_j$ is the denominator of $\frac{B_j}{4j}$, 
where $B_j$ is the $j^{\rm th}$ Bernoulli number.

In general, it is not known if the above result gives all restrictions
on the range of $L$. However, in special cases it does. For example,
for immersions $S^3\to\tR^5$ the range of $L$ is $\tZ$ (see
~\cite{E2} and also Remark ~\ref{lkrmk}) and for immersions
$S^{67}\to\tR^{69}$ it is $35\tZ$ (see Section ~\ref{rmksonl}).

In ~\cite{E1}, the author gives a complete classification of generic
immersions 
$M^{k}\to\tR^{2k-r}$, $r=0,1,2$, $k\ge 2r+4$ up to regular homotopy
through generic immersions (under some conditions on the lower
homotopy groups of $M$). 
Here, the class of a generic immersion is determined by its
self-intersection with induced natural additional  
structures (e.g. spin structures). The existence of invariants such as 
$L$ mentioned above implies that the corresponding classification in
other dimensions is more involved (see Remark ~\ref{nonloc}).

The invariants of generic immersions in Theorems
~\ref{thmL} and ~\ref{thmJ} give rise to 
invariants of regular homotopy which take values in finite cyclic
groups (Section ~\ref{rhinv}). We construct examples showing that,
depending on the source and target manifolds, the regular homotopy
invariants may or may not be trivial (Section ~\ref{expbm}). It would
be interesting to relate these invariants to invariants arising from
the cobordism theory of immersions (see for example Eccles ~\cite{Ec}).

\section{Statements of the main results}\label{mainres}
Before stating the main results we define the notions of invariants of
orders zero and one. They are similar to knot invariants of finite type
introduced by Vassiliev in ~\cite{V}.  

As in the Introduction, let $\tF$ denote the space of immersions and
let $\Si\subset\tF$ denote the discriminant hypersurface.
The set $\Si^1\subset\Si$ of non-generic immersions which appear at isolated
instances in generic 1-parameter families (see Lemmas 
~\ref{vdef1sm} and ~\ref{vdef1big} for descriptions of such immersions
$M^{nm-1}\to W^{n(m+1)-1}$) is a  
smooth submanifold of $\tF$ of codimension one. 
($\tF$ is an open subspace of the space of smooth maps, thus an
infinite dimensional manifold and the notion of codimension in $\tF$
makes sense.) 

If $f_0\in\Si^1$ then there is a neighborhood $U(f_0)$ of $f_0$ in
$\tF$ cut in two parts by 
$\Si^1$. A coherent choice of a positive and negative part of $U(f_0)$
for each $f_0\in\Si^1$ is a {\em coorientation} of $\Si$. (The
discriminant hypersurface in the space of immersions 
$M^{nm-1}\to W^{n(m+1)-1}$ will be cooriented in Section ~\ref{codiscr}.)   

Let $a$ be an invariant of generic immersions and let $f_0\in\Si^1$.
Define the jump $\nabla a$ of $a$ as
$$
\nabla a(f_0)=a(f_+)-a(f_-),
$$
where $f_+$ and $f_-$ are generic immersions in the positive
respectively negative part of $U(f_0)$. Then $\nabla a$ is a locally
constant function on $\Si^1$.
 
An invariant $a$  of generic immersions is a {\em zero order
invariant} if $\nabla a\equiv0$.  

Let $\Si^2$ be the set of all immersions which appear in generic
2-parameter families but can be avoided in generic 1-parameter
families (see Lemma ~\ref{vdef2} for the case 
$M^{nm-1}\to W^{n(m+1)-1}$). Then $\Si^2\subset\Si$ is a smooth
codimension two submanifold of $\tF$. 

An invariant $a$ of generic immersions is a {\em first order
invariant} if $\nabla a(f_0)=\nabla a(f_1)$ for any immersions
$f_0,f_1\in\Si^1$ which can be joined by a path in 
$\Si^1\cup\Si^2$ such that at intersections with $\Si^2$ its tangent
vector is transversal to the tangent space of $\Si^2$.

\begin{thm}\label{thmL}
Let $m>1$ and $n>1$ be integers and let $M^{nm-1}$ be a closed
manifold.  
\begin{itemize}
\item[{\rm (a)}] If $m$ is even and
$H_{n-1}(M;\tZ_2)=0=H_{n}(M;\tZ_2)$ then there exists a unique (up to 
addition of zero order invariants) first order $\tZ_2$-valued
invariant $\Lambda$ of generic 
immersions $M^{nm-1}\to \tR^{n(m+1)-1}$ satisfying the following
conditions:  It jumps by $1$ on the part of $\Si^1$ which consists
of immersions with one $(m+1)$-fold self-intersection point and  does
not jump on other parts of $\Si^1$. 
\item[{\rm (b)}] If $n$ is even, $M$ is orientable, and 
$H_{n-1}(M;\tZ)=0=H_{n}(M;\tZ)$ then there exists a unique (up to 
addition of zero order invariants) first order integer-valued
invariant $L$ of generic 
immersions $M^{nm-1}\to \tR^{n(m+1)-1}$ satisfying the following
conditions:  It jumps by $m+1$ on the part of $\Si^1$ which consists
of immersions with one $(m+1)$-fold self-intersection point and  does
not jump on other parts of $\Si^1$. 
\end{itemize}
\end{thm}

Theorem ~\ref{thmL} is proved in Section ~\ref{pfthmL}. The invariants 
$\Lambda$ and $L$ are defined in Definition ~\ref{dfnLambda} and
Definition ~\ref{dfnL}, respectively. If appropriately normalized,
$\Lambda$ and $L$ are additive under connected summation of generic
immersions (see Section ~\ref{csgi}).  

\begin{thm}\label{thmdivL}
Let $f\colon S^{2m-1}\to\tR^{2m+1}$ be a generic immersion. 
\begin{itemize}
\item[{\rm (a)}] If $m=4j+1$ then $2j+1$ divides $L(f)$.
\item[{\rm (b)}] If $m=4j+3$ then $4j+4$ divides $L(f)$.
\item[{\rm (c)}] If $m=2j$ then $p^r$  divides $L(f)$ for every prime
                 $p$ and integer $r$ such that $p^{r+k}$ divides
                 $2j+1$ and $p^{k+1}$ does not divide $\mu_j$ for 
                 some integer $k$. Here $\mu_j$ is the denominator of
                 $\frac{B_j}{4j}$, where $B_j$ is the $j^{\rm th}$
                 Bernoulli number. 
\end{itemize}
\end{thm}

We prove Theorem ~\ref{thmdivL} in Section ~\ref{pfthmdivL}. 

\begin{thm}\label{thmJ}
Let $m>1$ and $n>1$ be integers and let $M^{nm-1}$ be a closed
manifold and let $W^{n(m+1)-1}$ be a manifold. 
\begin{itemize}
\item[{\rm (a)}] If $n$ is odd then, for each integer $2\le r\le m$
such that $m-r$ is even, there exist a unique (up to addition of zero
order invariants) 
first order integer-valued invariant $J_r$,  of generic immersions 
$M^{nm-1}\to W^{n(m+1)-1}$ satisfying the following conditions:  
The invariant jumps by $2$ on the part of $\Si^1$ which consists
of immersions with one degenerate $r$-fold self-intersection point and
does not jump on other parts of $\Si^1$. 
\item[{\rm (b)}] If $n=2$ then there exists a unique (up to
addition of zero order invariants) first order integer-valued
invariant $J$,  of 
generic immersions $M^{2m-1}\to W^{2m+1}$ satisfying the following
conditions:   
The invariant jumps by $1$ on the part of $\Si^1$ which consists
of immersions with one degenerate $m$-fold self-intersection point and
does not jump on other parts of $\Si^1$. 
\end{itemize}
\end{thm}

Theorem ~\ref{thmJ} is proved in Section ~\ref{pfthmJ}. 
If appropriately normalized, $J_r$ and $J$ are additive under
connected summation of generic immersions (see Section ~\ref{csgi}).
The value of $J_r$
on a generic immersion is the Euler characteristic of its resolved
$r$-fold self-intersection manifold, the value of $J$ is the number of 
components in its $m$-fold self-intersection (which is a closed
1-dimensional manifold).

\section{Generic immersions and generic regular deformations}
In this section we define generic immersions and describe the
immersions corresponding to non-generic instances in generic 1- and
2-parameter families. We also describe the versal deformations of
these non-generic immersions. 

\subsection{Generic immersions}
Let $M$ and $W$ be smooth manifolds and let
$f\colon M\to W$ be an immersion. A point $q\in W$ is a {\em $k$-fold
self-intersection point of $f$} if $f^{-1}(q)$ consists of exactly $k$
points. Let  
$\Gamma_k(f)\subset W$ denote the set of $k$-fold self-intersection
points of $f$ and let $\wt {\Gamma}_k(f)=f^{-1}(\Gamma_k(f))\subset M$
denote its preimage. Note that $f|{\wt{\Gamma}_k(f)}\to\Gamma_k(f)$
is a $k$-fold covering. 

\begin{dfn}\label{dfngi}
Let $m>1$ and $n>1$. An immersion 
$f\colon M^{nm-1}\to W^{n(m+1)-1}$ is {\em generic} if 
\begin{itemize}
\item[{\bf g1}] $\Gamma_k(f)$ is empty for $k>m$, and
\item[{\bf g2}] if $q=f(p_1)=\dots=f(p_r)\in\Gamma_r(f)$ $(r\le m)$
then for any $i$, $1\le i\le r$ 
$$
df\, T_{p_i}M+\cap_{s\ne i}df\, T_{p_s}M=T_q W.
$$
\end{itemize}
\end{dfn}
A standard application of the jet-transversality theorem shows that
the set of generic immersions is open and dense in the space of all
immersions $\tF$. 

We remark that the set of $k$-fold self-intersection points
$\Gamma_k(f)$ of a generic immersion 
$f\colon M^{nm-1}\to W^{n(m+1)-1}$ is a smooth submanifold of
$W$ of dimension $n(m-k+1)-1$. Moreover, the deepest
self-intersection $\Gamma_m(f)$ is a closed manifold. 

Lemma ~\ref{giloc} below gives a 
local coordinate description of a generic immersion close to a
$k$-fold self-intersection point. To state it we
introduce some notation:  Let $x_i\in\tR^{nm-1}$, we
write $x_i=(x_i^0,x_i^1,\dots,x_i^{k-1})$, where
$x_i^0\in\tR^{n(m-k+1)-1}$ and $x_i^r\in\tR^{n}$, for $1\le r\le k-1$.
Similarly, we write $y\in\tR^{n(m+1)-1}$ as $y=(y^0,y^1,\dots,y^k)$, where
$y^0\in\tR^{n(m-k+1)-1}$ and $y^r\in\tR^{n}$, for $1\le r\le k$.
\begin{lma}\label{giloc}
Let $f\colon M^{nm-1}\to W^{n(m+1)-1}$ be a generic immersion and let
$q=f(p_1)=\dots=f(p_k)$ be a $k$-fold self intersection point of
$f$. Then 
there are coordinates $y$ in $V\subset W^{n(m+1)-1}$ centered at $q$
and coordinates $x_i$ on $U_i\subset M^{nm-1}$ centered at $p_i$,
$1\le i\le k$ such that $f$ is given by
\begin{align*}
f(x_1)&=(x_1^0,x_1^1,\dots,x_1^{k-2},x_1^{k-1},0),\\
f(x_2)&=(x_2^0,x_2^1,\dots,x_2^{k-2},0,x_1^{k-1}),\\
{}&\vdots{}\\
f(x_k)&=(x_k^0,0,x_k^1,x_k^2,\dots,x_k^{k-1}).
\end{align*}
(That is, if $y=f(x_i)$ then $y^0(x_i)=x_i^0$, $y^{r}(x_i)=x_i^r$ for
$1\le r\le k-i$, $y^{k-i+1}(x_i)=0$, and $y^{r}(x_i)=x_i^{r-1}$ for
$k-i+2\le r\le k$.)
\end{lma}
\begin{pf}
The proof is straightforward. 
\end{pf}

When the source and target of a generic immersion are oriented and the
codimension is even then there are induced orientations on the
self intersection manifolds: 
\begin{prp}
Let $n=2j>1$, $m>1$, $2\le k\le m$, and let 
$f\colon M^{2jm-1}\to W^{2j(m+1)-1}$ be a generic immersion. 
Orientations on $M^{2jm-1}$ and $W^{2j(m+1)-1}$ induce an orientation on
$\Gamma_k(f)$. 
\end{prp} 
\begin{pf}
Let $N$ denote the normal bundle of the immersion. The
decomposition
$$
f^*TW^{2j(m+1)-1}=TM^{2jm-1}\oplus N
$$ 
induces an orientation on $N$. 
If $q\in\Gamma_k(f)$, $q=f(p_1)=\dots=f(p_k)$ then 
$$
T_qW^{2j(m+1)-1}=T_q\Gamma_k(f)\oplus_{i=1}^k N_{p_i}. 
$$
The orientation on $N$ induces an orientation on $\oplus_{i=1}^k
N_{p_i}$. Since the dimension of the bundle $N$ is $2j$ which is even,
the orientation on the sum is independent on the ordering of the
summands. Hence, the decomposition above induces a well-defined
orientation on $T\Gamma_k(f)$.  
\end{pf}

\subsection{The codimension one part of the discriminant hypersurface}
Our next result describes (the 1-jets of) immersions in $\Si^1$ (see
Section ~\ref{mainres}).   
\begin{lma}\label{1jetSi1}
If $f_0\colon M^{nm-1}\to W^{n(m+1)-1}$ is an immersion in $\Si^1$ then 
{\bf g1} and {\bf g2} of Definition \ref{dfngi} holds, except at one
$k$-fold  
$(2\le k\le m+1)$ self-intersection point $q=f_0(p_1)=\dots=f_0(p_k)$,
where, 
\begin{itemize}
\item[{\rm (a)}] if $k=2$, 
$$
\dim(df_0\, T_{p_1}M+df_0\,
T_{p_2}M)=n(m+1)-2,
$$
or
\item[{\rm (b)}] if $2<k\le m+1$, for $i\ne l$ 
$$
\dim(df_0\, T_{p_i}M +\cap_{r\ne i,r\ne l}df_0\,
T_{p_r}M)=n(m+1)-1
$$
and 
$$
\dim(df_0\, T_{p_i}M+\cap_{r\ne i}df_0\,
T_{p_r}M)=n(m+1)-2.
$$
\end{itemize}
\end{lma} 
\begin{pf}
We have to show that degeneracies as above appears at
isolated parameter values in generic 1-parameter families, and that
further degeneracies (of the 1-jet) may be avoided . This follows
from the jet-transversality theorem applied to maps
$M^{nm-1}\times [-\delta,\delta]\to W^{n(m+1)-1}$. 
\end{pf}

A $k$-fold self-intersection point $q$ of an immersion 
$f_0\colon M^{nm-1}\to W^{n(m+1)-1}$ where {\bf g1} and {\bf g2} of
Definition ~\ref{dfngi} does not hold  will be called a  
{\em degenerate self-intersection point of $f_0$}. If $2\le k\le m$ we
say that $q$ is a  
{\em $k$-fold self-tangency point of $f_0$} if $k=m+1$ we say that $q$
is a {\em $(m+1)$-fold self-intersection point of $f_0$}. 

Recall that a deformation $F$ of a map $f$ is called {\em versal} if
any deformation of $f$ is equivalent (up to left-right action of
diffeomorphisms) to one induced from $F$. 

Let $f_0$ be an immersion in $\Si^1$. Then its versal deformation
$f_t$ is a 1-parameter 
deformation. In other words, it is a path $\lambda(t)=f_t$ in $\tF$
which intersects $\Si^1$ transversally at $f_0$ and thus,
$f_t$ are generic immersions for small $t\ne 0$.

Next, we shall describe immersions $f_0\in\Si^1$ in local coordinates
close to their degenerate self intersection point. 
Self-tangency points and $(m+1)$-fold self-intersection points are
treated in Lemma ~\ref{vdef1sm} and Lemma ~\ref{vdef1big},
respectively. 

To accomplish this we need a more detailed description of coordinates
than that given in Lemma ~\ref{giloc}, we use coordinates as there
with one more ingredient: 
We write (when necessary) $x_i^r\in\tR^n$ and $y_i^r\in\tR^n$, $r> 1$
as $x_i^r=(\xi_i^r,u_i^r)$ and 
$y^r=(\eta^r,v^r)$, where $\xi_i^r,\eta^r\in\tR$ and
$u_i^r,v_i^r\in\tR^{n-1}$. (Greek letters for scalars and Roman for
vectors). 

\begin{lma} \label{vdef1sm}
Let $f_0\colon M^{nm-1}\to W^{n(m+1)-1}$ be an immersion in $\Si^1$ 
and let $q=f_0(p_1)=\dots=f_0(p_k)$ be a point of degenerate $k$-fold
self intersection $2\le k \le m$. Then there are coordinates $y$ on
$V\subset W$
centered at $q$ and coordinates $x_i$ on $U_i\subset M$ centered at
$p_i$, $1\le i\le k$ such that in these coordinates the versal
deformation $f_t$, $-\delta<t<\delta$  of $f_0$ is constant 
outside of $\cup_i U_i$ and in neighborhoods of $p_i\in U_i$ 
it is given by 
\begin{itemize}
\item[{\rm (a)}] 
\begin{align*}
f_t(x_1)&=(x_1^0,x_1^1,\dots,x_1^{k-2},x_1^{k-1},0),\\
f_t(x_2)&=(x_2^0,x_2^1,\dots,x_2^{k-2},0,x_1^{k-1}),\\
{}&\vdots{}\\
f_t(x_{k-1})&=(x_{k-1}^0,x_{k-1}^1,0,x_{k-1}^2,\dots,x_{k-1}^{k-1}),
\end{align*}
(That is, if $y=f_t(x_i)$ then {\em for $1\le i\le k-1$},
$y^0(x_i)=x_i^0$, $y^1(x_i)=x_i^1$, $y^{r}(x_i)=x_i^r$ for
$1\le r\le k-i$, $y^{k-i+1}(x_i)=0$, and $y^{r}(x_i)=x_i^{r-1}$ for
$k-i+2\le r\le k$.)
\item[{\rm (b)}]
\begin{align*}
f_t(x_k)=&
\left(\begin{matrix}
x_k^0, & (\xi_k^1,0), & (\xi_k^2,u_k^1), & \dots, 
& (\xi_k^{k-1},u_k^{k-2}), & (-\xi_k^2-\dots-\xi_k^{k-1},u_k^{k-1}) 
\end{matrix}\right)\\
+&
\left(\begin{matrix}
0, & (0,0), & (0,0), &\dots, & (0,0), & (Q(x_k^0,\xi_k^1)+t,0)
\end{matrix}\right)
\end{align*}
where $Q$ is a nondegenerate quadratic form in the $n(m-k+1)$
variables $(x_k^0,\xi_k^1)$.
\end{itemize}
\end{lma}
\begin{pf}



It is straightforward to see that
we can find coordinates $x_i$ on $U_i$, $i=1,\dots,k$, so that up to
first order of approximation $f_0|U_i$ is given by the expressions in
(a) and the first term in (b) above.

We must consider also second order terms: Let $N$ be a linear subspace
in the coordinates $y$ transversal to $Tf_0(U_1)\cap\dots\cap
Tf_0(U_{k-1})$. Let $\phi\colon f(U_k)\to N$ be orthogonal projection
onto $N$. Then 
$\krn(d\phi)=Tf_0(U_1)\cap\dots\cap Tf_0(U_{k})$ 
and $\cokrn(d\phi)$ is 
1-dimensional. The second derivative
$d^2\phi\colon\krn(\phi)\to\cokrn(\phi)$    
must be a non-degenerate quadratic form, otherwise, we can avoid $f_0$ in
generic 1-parameter families. (This is a consequence of the jet
transversality theorem.) The Morse lemma then implies that, after
possibly adjusting the coordinates in $U_k$ by adding quadratic
expressions in $(x_k^0,\xi_k^1)$ to $\xi_k^j$, $j>1$, there exists
coordinates for $f_0$ as stated. 

Finally, we must prove that the deformation $f_t$ as given above is
versal. A standard result in singularity theory says that it is enough 
to prove that the deformation is infinitesimally versal. This is
straightforward.  
\end{pf}

In Lemma ~\ref{vdef1big} below we will use coordinates as at generic
$m$-fold self intersection  points. That is,
$x=(x^0,x^1,\dots,x^{m-1})\in\tR^{nm-1}$ and 
$y=(y^0,y^1,\dots, y^m)\in\tR^{n(m+1)-1}$, where $x^0,y^0\in\tR^{n-1}$ 
and $x_k,y_k\in\tR^n$ for $k\ne 0$.
\begin{lma}\label{vdef1big}
Let $f_0\colon M^{nm-1}\to W^{n(m+1)-1}$ be an immersion in $\Si^1$ 
and let $q=f_0(p_1)=\dots=f_0(p_{m+1})$ be a point of $(m+1)$-fold
self-intersection.  Then there are coordinates $y$ on $V\subset W$
centered at $q$ and coordinates $x_i$ on $U_i\subset M$ centered at
$p_i$, $1\le i\le k$ such that in these coordinates the versal
deformation $f_t$, $-\delta<t<\delta$  of $f_0$ is constant 
outside of $\cup_i U_i$ and in neighborhoods of $p_i\in U_i$ 
it is given by 
\begin{itemize}
\item[{\rm (a)}] 
\begin{align*}
f_t(x_1)&=(x_1^0,x_1^1,\dots,x_1^{m-2},x_1^{m-1},0),\\
f_t(x_2)&=(x_2^0,x_2^1,\dots,x_2^{m-2},0,x_1^{m-1}),\\
{}&\vdots{}\\
f_t(x_{m})&=(x_{m}^0,0,x_{m}^1,x_{m}^2,\dots,x_{m}^{m-1}),
\end{align*}
(That is, if $y=f_t(x_i)$ then {\em for $1\le i\le m$},
$y^0(x_i)=x_i^0$, $y^{r}(x_i)=x_i^r$ for
$1\le r\le m-i$, $y^{m-i+1}(x_i)=0$, and $y^{r}(x_i)=x_i^{r-1}$ for
$m-i+2\le r\le m$.)
\item[{\rm (b)}]
\begin{align*}
&f_t(x_k)=\\
&\left(\begin{matrix}
0, & (\xi_k^1,u^0_k), & (\xi_k^2,u_k^1), & \dots, 
& (\xi_k^{m-1},u_k^{m-2}), & (-\xi_k^1-\dots-\xi_k^{m-1}+t,u_k^{m-1}) 
\end{matrix}\right)
\end{align*}
\end{itemize}
\end{lma}

\begin{pf}
The proof is similar to the proof of Lemma ~\ref{vdef1sm}, but easier.
\end{pf}

\subsection{The codimension two part of the discriminant hypersurface}
\begin{lma}\label{vdef2}
Let $f_{0,0}\colon M^{nm-1}\to W^{n(m+1)-1}$ be an immersion in
$\Si^2$. Then either {\rm (a)} or {\rm (b)} below holds.
\begin{itemize} 
\item[{\rm (a)}] $f_{0,0}$ has two distinct degenerate
self-intersection points $q_1$ and $q_2$. Locally around $q_i$,
$i=1,2$ 
$f_{0,0}$ is as in Lemma ~\ref{vdef1sm} or Lemma ~\ref{vdef1big}.  
The versal deformation of $f_{0,0}$ is a product of the corresponding
1-parameter versal deformations.
\item[{\rm (b)}] $f_{0,0}$ has one degenerate $k$-fold $2\le k\le m$
self-intersection point $q$. There are coordinates $y$ centered at $q$ 
and coordinates $x_i$ centered at $p_i$, $1\le i\le k$ such that in
these coordinates the versal deformation $f_{s,t}$, $-\delta<t,s<\delta$  of
$f_{0,0}$ is constant outside of $\cup_i U_i$ and in a neighborhood of
$p_i\in U_i$ it is given by 
\begin{itemize}
\item[{\rm (b1)}] 
\begin{align*}
f_{s,t}(x_1)&=(x_1^0,x_1^1,\dots,x_1^{k-2},x_1^{k-1},0),\\
f_{s,t}(x_2)&=(x_2^0,x_2^1,\dots,x_2^{k-2},0,x_1^{k-1}),\\
{}&\vdots{}\\
f_{s,t}(x_{k-1})&=(x_{k-1}^0,x_{k-1}^1,0,x_{k-1}^2,\dots,x_{k-1}^{k-1}),
\end{align*}
That is, if $y=f_{s,t}(x_i)$ then {\em for $1\le i\le k-1$},
$y^0(x_i)=x_i^0$, $y^1(x_i)=x_i^1$, $y^{r}(x_i)=x_i^r$ for
$1\le r\le k-i$, $y^{k-i+1}(x_i)=0$, and $y^{r}(x_i)=x_i^{r-1}$ for
$k-i+2\le r\le k$.
\item[{\rm (b2)}]
\begin{align*}
&f_{s,t}(x_k)=\\
&\left(\begin{matrix}
x_k^0, & (\xi_k^1,0), & (\xi_k^2,u_k^1), & \dots, 
& (\xi_k^{k-1},u_k^{k-2}), & (-\xi_k^2-\dots-\xi_k^{k-1},u_k^{k-1}) 
\end{matrix}\right)\\
+&
\left(\begin{matrix}
0, & (0,0), & (0,0), &\dots, & (0,0), & (Q(x_k^0)+\xi_k^1((\xi_k^1)^2+s)+t,0)
\end{matrix}\right)
\end{align*}
\end{itemize}
where $Q$ is a nondegenerate quadratic form in the $n(m-k+1)-1$
variables $x_k^0$.
\end{itemize}
\end{lma}
\begin{pf}
The jet-transversality theorem applied to maps of 
$M^{2jm-1}\times [-\delta,\delta]^2$ into $W^{2j(m+1)-1}$ shows that
immersions with points of $k$-fold self-intersection $k\ge m+2$, as
well as immersions with points of $k$-fold self-intersection points
$2\le k\le m+1$ at which the 1-jet has further  
degenerations than the 1-jets of Lemma ~\ref{1jetSi1} can be avoided in
generic 2-parameter families. 

Assume that $f_{0,0}$ has a degenerate $k$-fold self-intersection
point. If $k=m+1$ it is easy to see that $f_{0,0}$ has the same local
form as the map in Lemma ~\ref{vdef1big}. So, immersions with $(m+1)$-fold
self-intersection points appears along 1-parameter subfamilies in
generic 2-parameter families. If $k<m+1$ we proceed as in the proof of
Lemma ~\ref{vdef1sm} and 
construct the map $\phi$. Applying the jet transversality theorem once 
again we see that if the rank of $d^2\phi$ is smaller than
$n(m-k+1)-1$ then $f_{0,0}$ can be avoided in generic 2-parameter
families. If the rank of $d^2\phi$ equals $n(m-k+1)-1$ then the third 
derivative of $\phi$ in the direction of the null-space of $d^2\phi$
must be non-zero otherwise $f_{0,0}$ can be avoided in generic
2-parameter families. 

With this information at hand we can find coordinates as
claimed. Finally, as in the proof of Lemma ~\ref{vdef1sm}, we must
check that the deformation given is infinitesimally versal. This is
straightforward. 
\end{pf}

\begin{figure}[htbp]
\begin{center}
\includegraphics[angle=0, width=6cm]{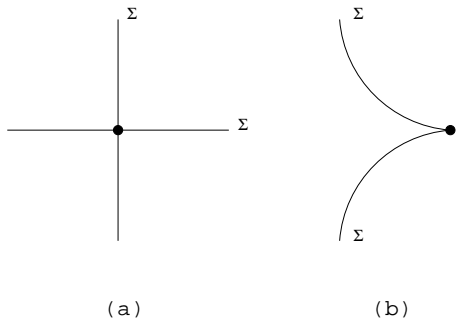}
\end{center}
\caption{The discriminant intersected with a small generic
2-disk.}\label{figcod2} 
\end{figure}
Pictures (a) and (b) in Figure \ref{figcod2} correspond to cases (a) and 
(b) in Theorem \ref{vdef2}. The codimension two parts of the
discriminant are represented by points. 
In (a) two branches of the discriminant, which consist of immersions
with one degenerate $r$-fold and one degenerate $k$-fold
self-intersection point respectively, intersect in $\Si^2$. In (b)
the smooth points of 
the semi-cubical cusp represents immersions with one degenerate
$k$-fold self-intersection point. The singular point represents
$\Si^2$. If an immersion in $\Si$ above the singular 
point ($\Si^2$) is moved below it then the index of the quadratic form
$Q$ in local coordinates close to the degenerate self-intersection
point (see Lemma \ref{vdef1sm}) changes.  

\section{Definition of the invariants}
In this section we define the invariants $J_r$, $J$, $L$ and
$\Lambda$. To this end, we describe resolutions of
self-intersections (for the definitions of $J_r$ and $J$) and
we compute homology of image-complements of and
of normal bundles (for the definitions of $L$ and $\Lambda$). 

\subsection{Resolution of the self intersection}
For generic immersions $f\colon M^{nm-1}\to W^{n(m+1)-1}$ let 
$
\Gamma^j(f)=\Gamma_j(f)\cup\Gamma_{j+1}(f)\cup\dots\cup\Gamma_m(f)
$, for $2\le j\le m$ and let $\wt{\Gamma}^j(f)=f^{-1}(\Gamma^j(f))$.
Resolving $\Gamma^j(f)$ we obtain a smooth manifold $\Delta^j(f)$: 
\begin{lma}\label{resol}
Let $m>1$ and $n>1$ and $2\le j\le m$ be integers. Let 
$f\colon M^{nm-1}\to W^{n(m+1)-1}$ be a generic immersion. Then there exists
closed $(n(m-j+1)-1)$-manifolds $\wt{\Delta}^j(f)$ and $\Delta^j(f)$,
unique up to diffeomorphisms, and immersions
$\sigma\colon\wt{\Delta}^j(f)\to M$ and $\tau\colon\Delta^j(f)\to W$
such that the diagram
$$
\begin{CD}
\wt{\Delta}^j(f) & @>{\sigma}>> & f^{-1}(\Gamma^j(f))\subset M\\
@VpVV &  & @VVfV\\
\Delta^j(f) & @>{\tau}>> & \Gamma^j(f)\subset W
\end{CD},
$$
commutes. The maps $\sigma$ and $\tau$ are surjective, have multiple
points only along $\wt{\Gamma}^{j+1}(f)$ and  $\Gamma^{j+1}(f)$
respectively, and $p$ is a $j$-fold cover.
\end{lma}
\begin{pf}
This is immediate from Lemma ~\ref{giloc}: Close to a $k$-fold self
intersection point $\Gamma_k(f)$ looks like the intersection of $k$
$(nm-1)$-planes in general position in $\tR^{n(m+1)-1}$. 
\end{pf}
\subsection{Definition of the invariants $J_r$ and $J$}
\begin{dfn}\label{dfnJ_r}
Let $m>1$ and $n>1$ be integers and assume that $n$ is odd. For a
generic immersion $f\colon M^{nm-1}\to W^{n(m+1)-1}$ and an integer
$2\le r\le m$ such that $m-r$ is even, define
$$
J_r(f)=\chi(\Delta^r(f)),
$$
where $\chi$ denotes Euler characteristic.  
\end{dfn}
\begin{dfn}\label{dfnJ}
Let $m>1$ be an integer. For a
generic immersion $f\colon M^{2m-1}\to W^{2m+1}$, define
$$
J_r(f)=\text{The number of components of $\Gamma_m(f)$.}
$$
\end{dfn}
\begin{lma}\label{wdJ}
The functions $J_r$ and $J$ are invariants of generic immersions.
\end{lma}
\begin{pf}
Let $f_t$, $0\le t\le 1$ be a regular homotopy through generic
immersions. If $F\colon M\times [0,1]\to W\times[0,1]$ is the
immersion $F(x,t)=(f_t(x),t)$ then
$\Gamma^j(F)\cong\Gamma^j(f_0)\times I$. It follows that
$\Delta^j(f_0)$ is diffeomorphic to $\Delta^j(f_1)$. 
\end{pf}

\subsection{Homology of complements of images}\label{imcom}
\begin{lma}\label{Hcom}
Let $f\colon M^{nm-1}\to\tR^{n(m+1)-1}$ be a generic
immersion. Assume that $M$ is closed. Then
\begin{itemize}
\item[{\rm (a)}]
$$
H_{n-1}(\tR^{n(m+1)-1}-f(M);\tZ_2)\cong H^{nm-1}(M;\tZ_2)\cong\tZ_2,
$$
and 
\item[{\rm (a)}] if $M$ is oriented then
$$
H_{n-1}(\tR^{n(m+1)-1}-f(M);\tZ)\cong H^{nm-1}(M;\tZ)\cong\tZ.
$$
\end{itemize}
\end{lma}
\begin{pf}
Alexander duality implies that 
$H_{n-1}(\tR^{n(m+1)-1}-f(M))\cong H^{nm-1}(f(M))$. By
Lemma \ref{giloc}, we can choose  
triangulations of $M$ and $f(M)$ such that the map $f\colon M\to f(M)$
is a bijection on the $(nm-k)$-skeleta for $1\le k\le n$.  
\end{pf}

\begin{rmk}\label{idHcom}
In case (a) of Lemma \ref{Hcom} we will use the unique
$\tZ_2$-orientation of  
$M$ and the duality isomorphism to identify
$H_{n-1}(\tR^{n(m+1)-1}-f(M);\tZ_2)$ with $\tZ_2$. In case (b) of
Lemma \ref{Hcom}, a $\tZ$-orientation of $M$  
determines an isomorphism $H^{nm-1}(M;\tZ)\to\tZ$. 
We shall use this and the duality isomorphism to
identify $H_{n-1}(\tR^{n(m+1)-1}-f(M);\tZ)$ with $\tZ$. 

A small ($n-1$)-dimensional sphere going around a fiber in the normal
bundle of $f(M)-\Gamma^2(f)$ generates these groups.  
\end{rmk}

\subsection{Homology of normal bundles of immersions}
Consider an immersion $f\colon M^{nm-1}\to W^{n(m+1)-1}$.
Let $N$ denote its normal bundle. Then $N$ is a 
vector bundle of dimension $n$ over $M$. Choose a Riemannian metric
on $N$ and consider the associated bundle $\partial N$ of unit vectors
in $N$. This is an ($n-1$)-sphere bundle over $M$. 
Let $\partial F\cong S^{n-1}$ denote a fiber of $\partial N$.

\begin{lma}\label{Hfib}
Let $f\colon M^{nm-1}\to W^{n(m+1)-1}$ be an immersion. Let
$i\colon \partial F\to\partial N$ denote the inclusion of the
fiber. 
\begin{itemize}
\item[{\rm (a)}] If $H_{n-1}(M;\tZ_2)=0=H_n(M;\tZ_2)$ then
$$
i_\ast\colon H_{n-1}(\partial F;\tZ_2)\to H_{n-1}(\partial N;\tZ_2),
$$
is an isomorphism. 
\item[{\rm (b)}] If $M$ and $W$ are oriented and
$H_{n-1}(M;\tZ)=0=H_n(M;\tZ)$ then
$$
i_\ast\colon H_{n-1}(\partial F;\tZ)\to H_{n-1}(\partial N;\tZ),
$$
is an isomorphism. 
\end{itemize}
\end{lma} 
\begin{pf}
This follows from the Leray-Serre spectral sequence.
\end{pf}
\begin{rmk}\label{idHfib}
In case (a) of Lemma \ref{Hfib} we use the isomorphism $i_\ast$ and a
canonical 
generator of $H_{n-1}(\partial F;\tZ_2)$ to identify 
$H_{n-1}(\partial N;\tZ_2)$ with $\tZ_2$. 
In case (b) of Lemma \ref{Hfib}, orientations of $M$ and $W$ induce an
orientation of the fiber sphere $\partial F$. We use this orientation
and $i_\ast$ in Lemma \ref{Hfib} to identify $H_{n-1}(\partial
N;\tZ)$ with $\tZ$.   
\end{rmk}

\subsection{Shifting the $m$-fold self-intersection}
Let $f\colon M^{nm-1}\to W^{n(m+1)-1}$ be a generic immersion.
Recall that the set of $m$-fold self-intersection points of $f$ is a
closed submanifold $\Gamma_m(f)$ of $W^{n(m+1)-1}$ of dimension
$n-1$ and its preimage $\wt{\Gamma}_m(f)$ is a closed submanifold of the 
same dimension in $M^{nm-1}$. 
\begin{lma}\label{sctn}
If $f\colon M^{nm-1}\to W^{n(m+1)-1}$ is an immersion then
there exists a smooth section  $s\colon\wt{\Gamma}_m(f)\to\partial N$,
where $\partial N$ is the unit sphere bundle of the normal bundle of
$f$. 
\end{lma}
\begin{pf}
This is immediate: The base space of the $n$-dimensional
vector bundle $N|\wt{\Gamma}_m(f)$ has dimension $n-1$.   
\end{pf} 
Let $\phi$ denote the canonical bundle map over
$f\colon M^{nm-1}\to W^{n(m+1)-1}$,
$$
\phi\colon N\to TW^{n(m+1)-1},
$$
and let $s$ be as in Lemma \ref{sctn}.
We define a vector field $v\colon \Gamma_m(f)\to TW^{n(m+1)-1}$
along $\Gamma_m(f)$: For $q=f(p_1)=\dots=f(p_m)$ let 
$$
v(q)=\phi(s(p_1))+\dots+\phi(s(p_m))).
$$
It is easy to see that $v$ is well-defined and smooth.

We now restrict attention to the case where the target manifold is
Euclidean space. 
\begin{dfn}\label{shif}
For a generic immersion $f\colon M^{nm-1}\to\tR^{n(m+1)-1}$,
let $\Gamma_m'(f,v,\epsilon)$ be the submanifold of $\tR^{n(m+1)-1}$
obtained by shifting $\Gamma_m(f)$ a small distance $\epsilon>0$ along $v$.
\end{dfn}
\begin{rmk}\label{skeps}
Note that for $\epsilon>0$ small enough $\Gamma_m'(f,v,\epsilon)\cap
f(M)=\emptyset$. Moreover, if $a>0$ is such that for all $0<\epsilon<a$ 
the corresponding $\Gamma_m'(f,v,\epsilon)$ is disjoint from $f(M)$ then
$\Gamma_m'(f,v,\epsilon)$ and $\Gamma_m'(f,v,a)$ are homotopic in
$\tR^{n(m+1)-1}-f(M)$. 
\end{rmk}
\subsection{Definition of the invariants $L$ and $\Lambda$}
\begin{dfn}\label{dfnLambda}
Let $f\colon M^{nm-1}\to\tR^{n(m+1)-1}$ be a generic immersion. 
Assume that $M$ satisfies 
\begin{equation*}
H_{n-1}(M;\tZ_2)=0=H_n(M;\tZ_2).
\tag{C$\Lambda$}
\end{equation*}
Define $\Lambda(f)\in\tZ_2$ as 
$$
\Lambda(f)=[\Gamma_m'(f,v,\epsilon)]-s_\ast[\wt{\Gamma}_m(f)]\in\tZ_2,
$$
where $\epsilon>0$ is very small and  
$
[\Gamma_m'(f,v,\epsilon)]\in H_{n-1}(\tR^{n(m+1)-1}-f(M);\tZ)\cong\tZ_2
$ 
(see Remark \ref{idHcom})
and 
$s_\ast[\wt{\Gamma}_m]\in H_{n-1}(\partial N;\tZ)\cong\tZ_2$ (see Remark
\ref{idHfib}). 
\end{dfn}
\begin{dfn}\label{dfnL}
Let $f\colon M^{nm-1}\to\tR^{n(m+1)-1}$ be a generic immersion. Assume
that $n$ is even and that $M$ is oriented and satisfies 
\begin{equation*}
H_{n-1}(M;\tZ)=0=H_n(M;\tZ).
\tag{C$L$}
\end{equation*}
Define 
$L(f)\in\tZ$ as  
$$
L(f)=[\Gamma_m'(f,v,\epsilon)]-s_\ast[\wt{\Gamma}_m(f)]\in\tZ,
$$
where $\epsilon>0$ is very small and  
$[\Gamma_m'(f,v,\epsilon)]\in
H_{n-1}(\tR^{2m(j+1)-1}-f(M);\tZ)\cong\tZ$  
(see Remark ~\ref{idHcom})
and 
$s_\ast[\wt{\Gamma}_m]\in H_{n-1}(\partial N;\tZ)\cong\tZ$ (see Remark
\ref{idHfib}). 
\end{dfn}
\begin{rmk}
If $\Tor(H_{n-2}(M;Z),\tZ_2)=0$
and
$M$ satisfies condition $({\rm C}L)$ in Definition ~\ref{dfnL}
 then by the universal coefficient
theorem $M$ also satisfies condition $({\rm C}\Lambda)$ in Definition
~\ref{dfnLambda}. It follows from Remarks ~\ref{idHcom},  \ref{idHfib} that
in this case $\Lambda=L\mod{2}$.  
\end{rmk}
\begin{lma}\label{wdL}
$\Lambda$ and $L$ are well-defined. That is, they do 
neither depend on the choice of $s$, nor on the choice of 
$\epsilon>0$.  
\end{lma}
\begin{rmk}
We shall often drop the awkward notation
$\Gamma_m'(f,v,\epsilon)$ and simply write $\Gamma'(f)$. This is
justified by Lemma ~\ref{wdL}.
\end{rmk}
\begin{pf}
The independence of $\epsilon>0$ follows immediately from Remark
~\ref{skeps}. We will therefore not write all $\epsilon$'s out in the
sequel of this proof.

Let $s_0$ and $s_1$ be two homotopic sections of 
$\partial N|\wt{\Gamma}_m(f)$. Let $v_0$ and $v_1$ be the
corresponding normal vector fields along $\Gamma_m(f)$. 

A homotopy $s_t$ between $s_0$ and $s_1$ induces a homotopy $v_t$ between
$v_0$ and $v_1$ and hence between $\Gamma_m'(f,v_0)$ and
$\Gamma_m'(f,v_1)$ in $\tR^{n(m+1)-1}-f(M)$. This shows that $\Lambda$
and $L$ only depend on the homotopy class of $s$.

To see that $\Lambda$ and $L$ are independent of the vector field we
introduce the notion of adding a local twist: Let 
$s\colon \wt{\Gamma}_m(f)\to\partial N$ be  
a section and let $p\in\wt{\Gamma}_m(f)$. Choose a neighborhood $U$ of $p$ in
$\wt{\Gamma}_m(f)$ and a trivialization $\partial N|U\cong U\times
S^{n-1}$ such that $s(u)\equiv w$, where $w$ is a point in
$S^{n-1}$. Let $D^{n-1}$ be disk inside $U$ and let $\sigma\colon D\to
S^{n-1}$ be a smooth map of degree $\pm 1$ such that $\sigma\equiv w$ in
a neighborhood of $\partial D$. Let $s^{\rm tw}$ be the section which is $s$
outside of $D$ and $\sigma$ in $D$. We say that $s^{\rm tw}$ is the result of
adding a local twist to $s$. (In the case when $M$ is
oriented and $n$ is even, the local twist is said to be positive if
the degree of 
$\sigma$ is $+1$ and negative if the degree of $\sigma$ is $-1$.)

Let $s^{\rm tw}$ be the vector field obtained by adding a twist to $s$. Let
$v^{\rm tw}$ and $v$ be the vector fields along $\Gamma_m(f)$ obtained
from $s^{\rm tw}$ and $s$ respectively. Clearly, 
$$
{s^{\rm tw}}_\ast[\Gamma_m(f)]=s_\ast[\Gamma_m(f)]\pm 1
\quad\text{ and }\quad 
[\Gamma_m'(f,v^{\rm tw})]=[\Gamma_m'(f,v)]\pm 1.
$$
Hence, $\Lambda$ and $L$ are invariant under adding local twists. 

A standard obstruction theory argument shows that 
if $s$ and $s'$ are two sections of $\partial N|\Gamma_m$ then by
adding local twists to $s'$ we can obtain a section $s''$ which is
homotopic to $s$. Hence, $\Lambda$ and $L$ are independent of the
choice of section. 
\end{pf}

\begin{lma}
$\Lambda$ and $L$ are invariants of generic immersions.
\end{lma}
\begin{pf}
Let $f_t$, $0\le t\le 1$ be a regular
homotopy through generic immersions. Consider the induced map 
$$
F\colon M\times I\to\tR^{n(m+1)-1}\times I,\quad F(x,t)=(f_t(x),t),
$$
where $I$ is the unit interval $[0,1]$. Shifting
$\Gamma_m(F)\cong\Gamma_m(f_0)\times I$  off $F(M\times I)$ using a
suitable 
vector field in $N|{\wt \Gamma}_m(F)$ it is easy to see that
$L(f_0)=L(f_1)$. 
\end{pf}

\section{Additivity properties, connected sum, and reversing orientation}
In this section we study how our 
invariants behave under two natural operations on generic immersions:
connected sum and reversing orientation.
\subsection{Connected summation of generic immersions}\label{csgi}
For (oriented) manifolds $M$ and $V$ of the same dimension, let
$M\skarp V$ denote the (oriented) connected sum of $M$ and $V$. 

Let $f\colon M^{nm-1}\to\tR^{n(m+1)-1}$ and $g\colon
V^{nm-1}\to\tR^{n(m+1)-1}$ be two generic immersions. 
We shall define the connected sum $f\skarp g$ of these. It will be a
generic immersion $M\skarp V\to\tR^{n(m+1)-1}$.

Let $(u,x)$,
$u\in\tR$ and $x\in\tR^{n(m+1)-2}$ be coordinates on
$\tR^{n(m+1)-1}$. Composing the immersions $f$ and $g$ with
translations we may assume that $f(M)\subset\{u\le-1\}$ and
$g(V)\subset\{u\ge 1\}$. Choose a point 
$p\in M$ and a point $q\in V$ such that there is only one point in
$f^{-1}(f(p))$ and in $g^{-1}(g(q))$. Pick an arc $\alpha$ in
$\tR^{n(m+1)-1}$ connecting
$f(p)$ to $g(q)$ and such that $\alpha\cap(f(M)\cup
g(V))=\{f(p),g(q)\}$. Moreover, assume that $\alpha$ meets $f(M)$ and
$g(V)$ transversally at its endpoints. Let $N$ be the normal bundle of 
$\alpha$. Pick a (oriented) basis of $T_{f(p)}f(M)$ and an
(anti-oriented) one 
of $T_{g(q)}(g(V))$. These give rise to $nm-1$ vectors over $\partial a$ 
in $N$. Extend these vectors to $nm-1$ independent normal vector
fields along $\alpha$. Using a suitable map of $N$ into a tubular
neighborhood of $\alpha$ these vector fields give rise to an embedding 
$\phi\colon\alpha\times D\to\tR^{n(m+1)-1}$, where $D$ denotes 
a disk of dimension $nm-1$, such that
$\phi|f(p)\times D$ is an (orientation preserving) embedding
into $f(M)$ and $\phi|g(q)\times D$ is an (orientation
reversing) embedding into $g(V)$. The tube 
$\phi(\alpha\times\partial D)$ now joins
$f(M)-\phi(f(p)\times\inr(D))$ to 
$g(V)-\phi(g(q)\times\inr(D))$. Smoothing the corners we get a 
generic immersion $f\skarp g\colon M\skarp V\to\tR^{n(m+1)-1}$.
\begin{lma}
Let $f,g\colon M^{nm-1}\to\tR^{n(m+1)-1}$ be generic immersions.
The connected sum $f\skarp g$ is independent
up to regular homotopy through generic immersions of both the choices
of points $f(p)$ and $g(q)$ and the choice of the path $\alpha$ used
to connect them.  
\end{lma}
\begin{pf}
This is straightforward.
(Note that the preimages of self-intersections has codimension $n>1$.)


\end{pf}

\begin{lma}
If $f,f'\colon M^{nm-1}\to\tR^{n(m+1)-1}$ and 
$g,g'\colon V^{nm-1}\to\tR^{n(m+1)-1}$, $m,n>1$ are regularly homotopic
through generic immersions then $f\skarp g$ and $f'\skarp g'$ are
regularly homotopic through generic immersions.
\end{lma}
\begin{pf}
This is straightforward
\end{pf}

\begin{prp}\label{csJ}
The invariants $J_r$ and $J$ are additive under connected summation.
\end{prp}
\begin{pf}
$\Delta^j(f\skarp g)=\Delta^j(f)\sqcup\Delta^j(g)$.
\end{pf}

Note that if $M^{nm-1}$ and $V^{nm-1}$ are manifolds which both
satisfy condition $({\rm C}\Lambda)$ in Definition ~\ref{dfnLambda} or
condition $({\rm C}L)$ in Definition ~\ref{dfnL} then so does $M\skarp V$. 

\begin{prp}\label{csL}
Let $f\colon M^{nm-1}\to\tR^{n(m+1)-1}$ and 
$g\colon V^{nm-1}\to\tR^{n(m+1)-1}$ be generic immersions. 
\begin{itemize}
\item[{\rm (a)}] If $M$ and $V$ both satisfy condition 
$({\rm C}\Lambda)$ then 
$$
\Lambda(f\skarp g)=\Lambda(f)+\Lambda(g).
$$
\item[{\rm (b)}] If $n$ is even and $M$ and $V$ are both oriented and
both satisfy condition $({\rm C}L)$ then 
$$
L(f\skarp g)=L(f)+L(g).
$$
\end{itemize}
\end{prp}
\begin{pf}
Note that $\Gamma_m(f\skarp g)=\Gamma_m(f)\sqcup\Gamma_m(g)$. 
Consider case (b). Choose
$2j$-chains $D$ and $E$ in $\{u<-1\}$ and $\{u>1\}$ bounding
$\Gamma_m(f)$ and $\Gamma_m(g)$ respectively and disjoint from the arc
$\alpha$, used in the construction of $f\skarp g$. Then 
$$
L(f\skarp g)=(D\cup E)\cdot f\skarp g(M\skarp V)=
D\cdot f(M)+E\cdot g(V)=L(f)+L(g).
$$ 
Case (a) is proved in exactly the same way.
\end{pf} 

\subsection{Changing orientation}
The invariants $J_r$, $J$, and $\Lambda$ are clearly orientation
independent. In contrast to this, the invariant $L$ is orientation
sensitive.  

To have $L$ defined, let $n=2j$ and consider an oriented closed manifold
$M^{2jm-1}$ which satisfies condition (C$L$).

\begin{prp}
Assume that there exists an orientation reversing diffeomorphism
$r\colon M\to M$. Let $f\colon
M^{2jm-1}\to\tR^{2j(m+1)-1}$ be a generic immersion. Then $f\circ r$
is a generic immersion and 
$$
L(f\circ r)=(-1)^{m+1}L(f)
$$
\end{prp}
\begin{pf}
Note that $\Gamma_m(f\circ r)=\Gamma_m(f)=\Gamma_m$. The orientation of
$\Gamma_m$ is induced from the decomposition
$$
T_q\tR^{2j(m+1)-1}=T_q\Gamma_m(f)\oplus N_{1}\oplus\dots\oplus N_{m}.
$$
The immersions $f$ and $f\circ r$ induces opposite orientations on
each $N_i$ hence the orientations induced on $\Gamma_m$ agrees if $m$
is even and does not agree if $m$ is odd. 
Let $D$ be a $2j$-chain bounding $\Gamma_m$, with its orientation
induced from $f$. If $m$ is even then
$$
L(f\circ r)=D\cdot f(r(M))=D\cdot -f(M)=-L(f).
$$ 
If $m$ is odd then
$$
L(f\circ r)=-D\cdot f(r(M))=-D\cdot -f(M)=L(f).
$$
\end{pf}
\begin{prp}
Let $R\colon \tR^{2j(m+1)-1}\to\tR^{2j(m+1)-1}$ be reflection in a
hyperplane. Let $f\colon M^{2jm-1}\to\tR^{2j(m+1)-1}$ be a generic
immersion. Then $R\circ f$ is a generic immersion and 
$$
L(R\circ f)=(-1)^mL(f)
$$
\end{prp}
\begin{pf}
Note that the oriented normal bundle of $Rf$ is $-RN$.


So the correctly oriented $2j$-chain bounding $\Gamma_{m}(R\circ f)$
is $(-1)^{m+1}RD$. Thus,
$$
L(R\circ f)=(-1)^{m+1} (RD\cdot Rf(M))=(-1)^m(D\cdot f(M))=(-1)^mL(f).
$$
\end{pf}
\section{Coorientations and proofs of Theorems ~\ref{thmL} and
~\ref{thmJ}} 
In this section we prove Theorems ~\ref{thmL} and ~\ref{thmJ}. To
do that we need a coorientation
of the discriminant hypersurface in the space of immersions.

\subsection{Coorienting the discriminant}\label{codiscr}
\begin{rmk}\label{nabla0}
Let $a$ be an invariant of generic immersions. If $a$ is
$\tZ_2$-valued then $\nabla a$ (see Section ~\ref{mainres}) is
well-defined without reference to any coorientation of
$\Si$. Moreover, if $a$ is integer-valued and $\Delta$ is a union of
path components of $\Si^1$  then the notion $\nabla a|\Delta\equiv0$
is well-defined without reference to a coorientation of $\Si$.
\end{rmk}

We shall coorient the relevant parts (see Remark ~\ref{nabla0}) of the
discriminant hypersurface. That is, we shall find coorientations of the
parts of the discriminant hypersurface in the space of generic
immersions where our invariants have non-zero jumps. These
coorientations will be {\em continuous} (see ~\cite{E3}, Section
7). That is, the intersection number of any generic small loop in
$\tF$ and $\Si^1$ vanishes and the coorientation extends continuously
over $\Si^2$. 

\begin{dfn}\label{coJr}
Let $m>1$ and $n>1$. Assume that $n$ is odd. Let $2\le
r\le m$ and assume that $m-r$ is even. Let $f_0\colon M^{nm-1}\to
W^{n(m+1)-1}$ be a generic immersion in $\Si^1$ with one degenerate
$r$-fold self-intersection point. Let $f_t$ be a versal deformation of 
$f_0$. We say that $f_{\delta}$ is on the positive side of $\Si^1$ and 
$f_{-\delta}$ on the negative side if
$\chi(\Delta^r(f_\delta))>\chi(\Delta^r(f_{-\delta}))$. 
\end{dfn}
\begin{rmk}
Note that by Lemma ~\ref{vdef1sm} $\chi(\Delta^r(f_\delta))$ is
obtained from  
$\chi(\Delta^r(f_{-\delta}))$ by a Morse modification. Since the
dimension of $\chi(\Delta^r(f_\delta))$ is $n(m-r+1)-1$ is even the
Euler characteristic changes under such modifications. 
\end{rmk}
\begin{dfn}\label{coJ}
Let $m>1$. Let $f_0\colon M^{2m-1}\to
W^{2m+1}$ be a generic immersion in $\Si^1$ with one degenerate
$m$-fold self-intersection point. Let $f_t$ be a versal deformation of 
$f_0$. We say that $f_{\delta}$ is on the positive side of $\Si^1$ and 
$f_{-\delta}$ on the negative side if the number of components in
$\Gamma_m(f_\delta)$ is larger than the number of components in
$\Gamma_m(f_{-\delta})$. 
\end{dfn}
\begin{rmk}
Note that by Lemma ~\ref{vdef1sm} $\Gamma_m(f_\delta)$ is
obtained from  
$\Gamma_m(f_{-\delta})$ by a Morse modification. Since these are
1-manifolds the number of components change.
\end{rmk}

The construction of the relevant coorientation for $L$ is more
involved. It is based on a high-dimensional counterpart of the   
notion of over- and under-crossings in classical knot theory. 

Let $f_0\colon M^{2j(m+1)-1}\to W^{2j(m+1)-1}$ be an immersion of
oriented manifolds. Assume that $f_0$ has an ($m+1$)-fold
self-intersection point,  
$q=f_0(p_1)=\dots=f_0(p_{m+1})$. Let $f_t$ be a versal deformation of
$f_0$. Let $U_i$ be small neighborhoods of $p_i$ 
and let $S_i^t$ denote the oriented sheet $f_t(U_i)$. Note that
$S_i^0\cap\dots\cap S^0_{m+1}=q$ and that 
$S^t_1\cap\dots\cap S^t_{m+1}=\emptyset$ if $t\ne 0$. 

Let $D^t_i=\cap_{j\ne i}S_j^t$. Let $w_i$ be a line transversal to
$TS_i^0+TD^0_i$ at $q$. For small $t\ne 0$ both $S_i^t$ and
$D_i^t$ intersects $w_i$. Orienting the line from $S_i^t$ to
$D_i^t$ gives a local orientation $(D_i^t,S_i^t,\vec{w}_i)$
of $\tR^{2j(m+1)-1}$. Comparing it with the standard orientation of
$\tR^{2j(m+1)-1}$ we get a sign
$\sigma_t(i)=\Or(D_i^t,S_i^t,\vec{w}_i)$, where $\Or$ denotes the
sign of the orientation. Note that, if we orient
the line from $D_i^t$ to $S_i^t$ we get the opposite orientation
$-\vec{w}_i$ of $w$ and
$$
\Or(D_i^t,S_i^t,\vec{w}_i)=\Or(S_i^t,D_i^t,-\vec{w}_i),
$$
since $D_i^t$ and $S_i^t$ are both odd-dimensional. Note also that
$\sigma_t(i)=-\sigma_{-t}(i)$ (see Lemma \ref{vdef1big}). 

We next demonstrate that $\sigma_t(i)=\sigma_t(j)$ for all $i,j$: 
Let $N_i^t$ be the oriented normal bundle of $S_i^t$. Let $w$ be a vector
from $D_1^t$ to $D_2^t$, transversal to both
$TS_1^0+TD^0_1$  and $TS_2^0+TD^0_2$. Orient it from
$D_1^t$ to $D_2^t$. Then 
$$
\sigma_t(1)=\Or(S_1^t,D^t_1,\vec w),
$$
and hence $TD^t_1+w$ gives a normal bundle $N_1^t$ and by the
convention used to orient the normal bundle
$\Or(N_1^t)=\sigma_t(1)\Or(D_1^t,\vec w)$. In a similar way it follows that
$\Or(N_2^t)=\sigma_t(2)\Or(D_2^t,-\vec w)$. 

Now, by definition
\begin{align*}
1= & \Or(D_1^t,N_2^t,\dots,N_m^t) 
=\sigma_t(2)\Or(D_1^t,(D_2^t,-\vec w),N_3^t,\dots,
N_m^t)=\\
= & \text{[$\dim(D_2^t)$ is odd]}= 
-\sigma_t(2)\Or(D_1^t,-\vec w,D_2^t,N_3^t,\dots, N_m^t)=\\
= & \sigma_t(2)\sigma_t(1)\Or(N_1^t,D_2^t,N_3^t,\dots,N_m^t)
=\sigma_t(1)\sigma_t(2).
\end{align*}
Hence, $\sigma_t(1)\sigma_t(2)=1$ as claimed.
\begin{dfn}\label{coL}
Let $f_0\colon M^{2j(m+1)-1}\to W^{2j(m+1)-1}$ be an immersion of
oriented manifolds. Assume that $f_0$ has an $(m+1)$-fold
self-intersection point. Let $f_t$ be a versal deformation of
$f_0$. We say  
that $f_{\delta}$ is on the positive side of $\Si^1$ at $f_0$ if 
$$
\sigma_\delta(1)=\dots=\sigma_\delta(m+1)=+1.
$$
\end{dfn}
\subsection{First order invariants}
The following obvious lemma will be used below.
\begin{lma}\label{jump}
Any first order invariant of generic immersions is uniquely (up to
addition of zero order invariants) determined by its jump.\qed 
\end{lma}
\subsection{Proof of Theorem ~\ref{thmJ}}\label{pfthmJ}
We know from Lemma ~\ref{wdJ} that $J$ and  $J_r$ are invariants of
generic immersions. We must calculate their jumps. 

We start in case (a): Assume that $n$ is odd.
Let $f_t\colon M^{nm-1}\to\tR^{n(m+1)-1}$, $t\in [-\delta,\delta]$ be
a versal deformation of $f_0\in\Si^1$ and fix $r$, $2\le r\le m$ such
that $m-r$ is even. 

If
$f_0$ has a degenerate $k$-fold intersection point $2\le k\le m+1$ and
$k\ne r$ then $\Delta^r(f_{-\delta})$ is diffeomorphic to
$\Delta^r(f_{\delta})$: If $k<r$ then $\Delta^r(f_{-\delta})$ is
not affected at all by the versal deformation. If $k>r$ then the
immersed submanifold $\tau^{-1}(\Gamma_k(f_{-\delta}))$ of
$\Delta^r(f_{-\delta})$ is changed by surgery (or is
deformed by regular homotopy if $k=m+1$) under the
versal deformation. This does not affect the diffeomorphism class of
$\Delta^r(f_{-\delta})$.      

If $f_0$ has a degenerate $r$-fold self-intersection point then
$\Delta^r{f_{\delta}}$ is obtained from $\Delta^r{f_{-\delta}}$ by a
surgery (see Lemma ~\ref{vdef1sm}). Since the dimension of
$\Delta^r(f_{-\delta})$ is $n(m-r+1)-1$ which is even this changes the 
Euler characteristic by $\pm2$.

According to our coorientation conventions $\nabla J_r(f_0)=2$ if
$f_0$ has a degenerate $r$-fold self-intersection point. Also, 
$\nabla J_r(f_0)=0$ if $f_0$ has any other degeneracy.

The jump of $J$ in case (b) is $1$ by the same argument. 

It is evident from Lemma ~\ref{vdef2} (see Figure
\ref{figcod2}) that $J$ and $J^r$ are first order invariants. The
theorem now follows  from Lemma ~\ref{jump}.\qed
\subsection{Proof of Theorem \ref{thmL}}\label{pfthmL}
We start with (b):
Let $f_t\colon M^{2jm-1}\to\tR^{2j(m+1)-1}$, $t\in [-\delta,\delta]$
be a generic one-parameter family intersecting $\Si^1$ in $f_0$. If
the degenerate intersection point of $f_0$ is a $k$-fold intersection
point with $2\le k\le m-1$ then clearly
$L(f_{-\delta})=L(f_{\delta})$, since the $m$-fold self-intersection
is not affected under such deformations.

Assume that $f_0$ has a degenerate $m$-fold self-intersection
point. Without loss of generality we may assume that $f_t$ is a 
deformation of the form in Lemma \ref{vdef1sm} (we use coordinates as
there): 

Let $F(x,t)=(f_t(x),t)$. Shifting $((\{Q(y^0,v^1)=t\},0\dots,0),t)$ we 
obtain a $2j$-chain bounded by 
$\Gamma_m'(f_\delta)\times\delta-\Gamma_m'(f_{-\delta})\times {-\delta}$ 
in
$\tR^{2m(j+1)-1}\times[-\delta,\delta]-F(M\times[-\delta,\delta])$. 
It follows that $L(f_\delta)=L(f_{-\delta})$.

If $f_0$ has an $(m+1)$-fold self-intersection then
the discussion preceding Definition ~\ref{coL} shows that
$L(f_{\delta})=L(f_{-\delta})+(m+1)$: At an $(m+1)$-fold self
intersection point $m+1$ crossings are turned into crossings of 
opposite sign.

Hence, $\nabla L(f_0)=m+1$ if $f_0$ has an $(m+1)$-fold
self-intersection point and $\nabla L(f_0)=0$ if $f_0$ has any other
degeneracy. 

The calculation of $\nabla \Lambda$ in (a) is analogous.  Let us
just make a remark about the parity of $m$: At an $(m+1)$-fold
self-intersection point $m+1$ crossings are  
changed. Hence, at instances of $(m+1)$-fold self-intersection the
invariant $\Lambda$ changes by $1\in\tZ_2$ if $m+1$ is odd and does not
change if $m+1$ is even. 

It is
immediate from Lemma ~\ref{vdef2}, see Figure \ref{figcod2} that
$\Lambda$ and $L$ are first order invariants. 
The theorem now follows from Lemma ~\ref{jump}.
\qed

\section{Invariants of regular homotopy}\label{rhinv}
A function of immersions $M\to W$ which is constant on path components 
of $\tF$ will be called an {\em invariant of regular homotopy}. Our
geometrically defined invariants of generic immersion give rise to
torsion invariants of regular homotopy.
\begin{dfn}
Let $n$ be odd. Let $f\colon M^{nm-1}\to W^{n(m+1)-1}$ be an
immersion. Let $f'$ be a generic immersion regularly homotopic to $f$.
For $2\le r\le m$ such that $m-r$ is even, define $j_r(f)\in\tZ_2$ as
$$
j_r(f)=J_r(f')\mod{2}.
$$
\end{dfn}
\begin{prp}\label{rhJ}
The function $j_r$ is an invariant of regular homotopy.
\end{prp}
\begin{pf}
Clearly it is enough to show that for any two regularly homotopic
generic 
immersions $f_0$ and $f_1$, $j_r(f_0)=j_r(f_1)$. Let $f_t$ be a
generic regular homotopy from $f_0$ to $f_1$. Then $f_t$ intersects
$\Si$ transversally in a finite number of points in $\Si^1$. It
follows from Theorem ~\ref{thmJ} that $j_r$ remains unchanged at such
intersections. Hence, $j_r(f_0)=j_r(f_1)$.
\end{pf}
\begin{dfn}
Let $n$ be even. Assume that $M^{nm-1}$ is a manifold which satisfy
condition $({\rm C}L)$. Let $f\colon M^{nm-1}\to \tR^{n(m+1)-1}$ be an
immersion. Let $f'$ be a generic immersion regularly homotopic to $f$.
Define $l(f)\in\tZ_{m+1}$ as
$$
l(f)=L(f')\mod{(m+1)}.
$$
\end{dfn}
\begin{prp}\label{rhL}
The function $l$ is an invariant of regular homotopy.
\end{prp}
\begin{pf}
The proof is identical to the proof of Proposition ~\ref{rhJ}.
\end{pf}
\begin{dfn}
Let $m>1$ be odd. Assume that $M^{nm-1}$ is a manifold which satisfy
condition $({\rm C}\Lambda)$. Let $f\colon M^{nm-1}\to \tR^{n(m+1)-1}$ be an
immersion. Let $f'$ be a generic immersion regularly homotopic to $f$.
Define $\lambda(f)\in\tZ_2$ as
$$
\lambda(f)=\Lambda(f').
$$
\end{dfn}
\begin{prp}
The function $\lambda$ is an invariant of regular homotopy.
\end{prp}
\begin{pf}
This follows from the proof of Theorem ~\ref{thmL}, where it is noted that
when $m$ is odd $\Lambda$ remains constant when a regular homotopy
intersects $\Si^1$. 
\end{pf}

\section{Sphere-immersions in codimension two}
In this section Theorem \ref{thmdivL} will be proved. To do that we 
will first discuss the classifications of sphere-immersions and
sphere-embeddings up to regular homotopy. 

\subsection{The Smale invariant}
In Smale's classical work \cite{S} it is proved that there is a
bijection between the set of regular homotopy classes $\Imm(k,n)$ of
immersions $S^k\to\tR^{k+n}$ and the
elements of the group $\pi_k(V_{{k+n},k})$, the $k^{\rm th}$ homotopy
group of the Stiefel manifold of $k$-frames in $(k+n)$-space. 
If $f\colon S^k\to\tR^{k+n}$ is an immersion we let
$\Omega(f)\in\pi_k(V_{{k+n},k})$ denote its Smale invariant.
Via $\Omega$ we can view $\Imm(k,n)$ as an Abelian group. 

If the codimension is two the groups appearing in Smale's
classification are easily computed: The exact homotopy
sequence of the fibration 
$$
SO(2)\hookrightarrow SO(2m+1)\to V_{2m+1,2m-1}
$$ 
implies that
$\pi_{2m-1}(V_{2m+1,2m-1})\cong\pi_{2m-1}(SO(2m+1))$. Bott-periodicity 
then gives:
$$
\pi_{2m-1}(V_{2m+1,2m-1})=\begin{cases}
\tZ &\text{if } m=2j,\\
\tZ_2 &\text{if } m=4j+1,\\
0 &\text{if } m=4j+3.
                          \end{cases}
$$
\begin{rmk}\label{geomope}
It is possible to identify the group operations in $\Imm(k,2)$ geometrically.
Kervaire \cite{K} proves that the Smale invariant $\Omega$ is
additive under connected sum of immersions. This gives the geometric
counterpart of addition. 
If $f\colon S^n\to\tR^{n+2}$ is an immersion, $r\colon S^n\to S^n$
is an orientation reversing diffeomorphism, and $n\ne 2$ then 
$\Omega(f)=-\Omega(f\circ r)$. (See ~\cite{E2} for the case
$S^3\to\tR^5$, the other cases are analogous). This gives
the geometric counterpart of the inverse operation in $\Imm(n,2)$ for
$n\ne 2$.  

It is interesting to note that for immersions $f\colon S^2\to\tR^4$, 
$\Omega(f)=\Omega(f\circ r)$: The Smale invariant can in this 
case be computed as the algebraic number of self-intersection points. 
This number is clearly invariant under reversing orientation. (The same is
true also for immersions $S^{2k}\to\tR^{4k}$, $k\ge 1$.) 
\end{rmk}
\begin{lma}\label{lhom}
The regular homotopy invariant $l$ induces a homomorphism
$$
\Imm(2m-1,2)\to\tZ_{m+1}.
$$
\end{lma}
\begin{pf}
Note that the dimensions are such that $L$ is defined and 
spheres certainly satisfy condition $({\rm C}L)$. The invariant $L$ is
additive under 
connected summation and changes sign if an immersion is composed on
the left with an orientation reversing diffeomorphism. Hence, $l$
induces a homomorphism. 
\end{pf}
Let $\Emb(n,2)\subset\Imm(n,2)$ be the set of regular homotopy classes
which contain embeddings. By Remark ~\ref{geomope}
$\Emb(n,2)$ is a subgroup of $\Imm(n,2)$. A result of Hughes and
Melvin ~\cite{HM} states that 
$\Emb(4j-1,2)\subset\Imm(4j-1,2)\cong\tZ$ is a subgroup of index
$\mu_j$. Here $\frac{B_j}{4j}=\frac{\nu_j}{\mu_j}$, where $B_j$ is the 
$j^{\rm th}$ Bernoulli number and $\nu_j$ and $\mu_j$ are coprime
integers.   

\subsection{Proof of Theorem ~\ref{thmdivL}}\label{pfthmdivL}
By Lemma \ref{lhom},
$l\colon\Imm(2m-1,2)\to\tZ/(m+1)\tZ$ is a homomorphism and clearly,
$\Emb(2m-1,2)\subset\krn(l)$. 

In case (b) $\Imm(8j+5,2)=0$ and hence the image of $l$ is zero, 
which proves that $L(f)$ is always divisible by $m+1=4j+4$. 

In case (a) $\Imm(8j+1,2)\cong\tZ/2\tZ$. Hence, for any immersion $f$,
$l(f\skarp f)=0$. Thus, $L(f\skarp f)=2L(f)$ is divisible by
$m+1=4j+2$, which implies that $L(f)$ is divisible by $2j+1$.

In case (c) $\Emb(4j-1,2)\cong \mu_j\tZ$ as a subgroup of
$\Imm(4j-1,2)\cong\tZ$. Hence, for any immersion $f$, $m+1=2j+1$
divides 
$$
L(f\skarp\dots[\mu_j \text{ summands}]\dots\skarp f)=\mu_jL(f).
$$ 
Thus, if  if $p$ is a prime and $r,k$ are integers such that $p^{r+k}$
divides $2j+1$ and $p^{k+1}$ does not divide $\mu_j$ then $p^{r}$
divides $L(f)$.
\qed

\section{Examples and problems}\label{expbm}
In this section we construct examples of 1-parameter families of
immersions which shows that all the first order invariants we have
defined are non-trivial. We also discuss some problems in connection
with the regular homotopy invariants defined in Section ~\ref{rhinv}.

\subsection{Examples}
Using our local coordinate description of immersions with one
degenerate self-intersection point we can construct examples showing
that the invariants $J$, $J_r$, $\Lambda$ and $L$ are non-trivial. We
start with $\Lambda$ and $L$:

Choose $m$ standard spheres $S_1,\dots,S_m$ of dimension $(nm-1)$
intersecting in general position in $\tR^{n(m+1)-1}$ so that
$S_1\cap\dots\cap S_m\cong S^{n-1}$. Pick a point 
$p\in S_1\cap\dots\cap S_m$ and let 
$S_{m+1}$ be another standard $(nm-1)$-sphere intersecting 
$S_1\cap\dots\cap S_m$ at $p$ so that in a neighborhood of $p$ the
embeddings are given by the expressions in Lemma ~\ref{vdef1big} and
so that $p$ is the only degenerate intersection point. Let
$f_0\colon S^{nm-1}\to\tR^{n(m+1)-1}$ be the immersions which is the
connected sum of $S_1,\dots,S_{m+1}$. Let $f_t$ be a versal
deformation of $f$. Then, after possibly reversing the direction of
the versal deformation we have $\Lambda(f_{\delta})=1\in\tZ_2$, if $n$ 
is odd or $L(f_{\delta})=\pm(m+1)$ if $n$ is even. 

In the latter case we would get the other sign of $L(f_{\delta})$ if
the 
orientation of $S_{m+1}$ in the construction above is reversed. To
distinguish these two immersions denote one by $h^+$ and the other by 
$h^-$, so that $L(h^+)=m+1=-L(h^-)$. Then it is an immediate
consequence of Theorem ~\ref{thmL} that: 

{\em Any generic regular homotopy from $h^+$ to $h^-$ has at least two 
instances of $(m+1)$-fold self-intersection.}

Theorem ~\ref{thmL} has many corollaries as the one just mentioned.
 
To see that $L$ and $\Lambda$ are nontrivial for other source
manifolds. We use connected sum with the immersions
$S^{nm-1}\to\tR^{n(m+1)-1}$ just constructed and Lemma ~\ref{csL}.

The proof that the invariants $J$ and $J_r$ are nontrivial are along
the same lines: Construct a sphere-immersion into Euclidean space with 
one degenerate self- intersection point as in the local models in Lemma
~\ref{vdef1sm}. Then embed this Euclidean space with immersed sphere 
into any target manifold and use connected sum together with Lemma
~\ref{csJ}. 
\begin{rmk}\label{Lrange}
Applying connected sum to the sphere-immersion constructed above in
the case when it is an immersion $S^{2m-1}\to\tR^{2m+1}$ shows that
for any integer $k$ there exists a  generic immersions 
$f\colon S^{2m-1}\to\tR^{2m+1}$ such that $L(f)=(m+1)k$.
\end{rmk}
\begin{rmk}\label{nonloc}
Let $f_0\colon M^{nm-1}\to\tR^{n(m+1)-1}$ be an immersion in $\Si^1$
with an $(m+1)$-fold self-intersection point and let $f_t$, $-\delta<
t<\delta$ be a versal deformation of $f_0$. 

Let $\Gamma_{\pm}=\cup_{i=1}^m\Gamma_{i}(f_{\pm\delta})$ and 
${\wt \Gamma}_{\pm}=f_{\pm\delta}^{-1}(\Gamma_{\pm})$. Then if
$U_{\pm}$ are small regular neighborhoods of $\Gamma_{\pm}$ and 
${\wt U}_{\pm}=f_{\pm\delta}^{-1}(U_{\pm})$ then there exist
diffeomorphisms $\phi$ and $\psi$ such that the following diagram
commutes
$$
\begin{CD}
({\wt U}_-,{\wt \Gamma}_-)  & @>\psi>> & ({\wt U}_+,{\wt \Gamma}_+)\\
@V{f_{-\delta}}VV & {} & @VV{f_{+\delta}}V\\
(U_-,\Gamma_-)  & @>\phi>> & (U_+,\Gamma_+)
\end{CD}.
$$
That is, the local properties of a generic immersion close to its 
self-intersection does not change at instances of $(m+1)$-fold
self-intersection points in generic 1-parameter families.

Assume that $M$ fulfills the requirements of Theorem \ref{thmL}.
Then $\Lambda$ or $L$ is defined and 
$\Lambda(f_{-\delta})\ne \Lambda(f_{\delta})$ or
$L(f_{-\delta})\ne L(f_{\delta})$, respectively. This implies that
$f_{-\delta}$ and $f_{\delta}$ are not regularly homotopic through
generic immersions. Hence, knowledge of the local properties of a
generic immersion $M^{nm-1}\to\tR^{n(m+1)-1}$ close to its
self-intersection is not enough to 
determine it up to regular homotopy through generic immersions.
\end{rmk}
\subsection{Problems}\label{pbm}

{\em Are the regular homotopy invariants $l$, $\lambda$ and $j_r$
nontrivial?} 

A negative answer to this question implies restrictions on
self-intersection manifold. A positive answer gives non-trivial
geometrically defined regular homotopy invariants. 

\subsection{Remarks on the invariant $l$}\label{rmksonl}
Theorem \ref{thmdivL} gives information on the possible range of $l$
for sphere-immersions in codimension two. We consider the most
interesting cases of immersions
$S^{4j-1}\to\tR^{4j+1}$, $j\ge 1$ (Theorem \ref{thmdivL} (c)).
In the first case $S^3\to\tR^5$, $l$ is non-trivial. This was shown in 
\cite{E2}. Hence, for any integer $b$ there are generic immersions
$f\colon S^3\to\tR^5$ such that $L(f)=b$. 
\begin{rmk}\label{lkrmk}
The invariant $l$ is called $\lambda$ in \cite{E2} and is the
$\mod{3}$ reduction of an integer-valued invariant called $\lk$.
The definition of $\lk$ given in ~\cite{E2} differs slightly from the 
definition of $L$ given here. There is an easy indirect way to see that,
nonetheless, the two invariants are the same: 
Due to Theorem ~\ref{thmL} and Theorem 2
in \cite{E2}, $L-\lk$ is an invariant of regular homotopy which is $0$
on embeddings and additive under connected summation. Assume that
$L(f)-L'(f)=a$ for some immersion. Then, since $\mu_1=24$, the
connected sum of $24$ copies of any immersion is regularly homotopic
to an embedding. Hence, $24a=0$ and therefore $a=0$. Thus, $L\equiv\lk$. 
\end{rmk}
If $2j+1$ is prime then Theorem \ref{thmdivL} does not impose any
restrictions on $l$ since, in this case, $2j+1$ divides $\mu_j$ (see
Milnor and Stasheff \cite{MS}).

On the other hand, if none of the prime factors of $2j+1$ divides
$\mu_j$, Theorem \ref{thmdivL} implies that $l$ is 
trivial and in such cases Remark \ref{Lrange} allows us to determine
the range of $L$ which is $(2j+1)\tZ$. The first two cases where this
happens are $S^{67}\to\tR^{69}$ and $S^{107}\to\tR^{109}$, where $2j+1$ equals
$35=5\cdot 7$ and $55=5\cdot 11$ and $\mu_j$ equals $24=2^3\cdot 3$ and  
$86184=2^3\cdot 3^4\cdot 7\cdot 19$, respectively.


\subsection{Remarks on the invariant $j_2$}
The first case in which to consider the invariant $j_2$ are 
immersions of a $5$-manifold into an $8$-manifold.

Theorem 7.30 in ~\cite{E1} shows that  
the self-intersection surface of a  generic immersions of
$S^5\to\tR^8$ must have {\em even} Euler characteristic (even though it may
be non-orientable). Thus, in this case $j_2$ is trivial. (Note that
$\pi_5(V_{8,5})\cong\tZ_2$ so that there are two regular homotopy
classes of immersions $S^5\to\tR^8$. Hence, it is non-trivial to see
that $l$ is trivial in this case.)   

In contrast, the invariant $j_2$ is non-trivial for immersions
$S^5\to\tR P^8$: Clearly, $S^5$ embeds in $\tR P^8$. Thus, there is
an immersion $f$ with $j_2(f)=0$.

Consider three hyperplanes $H_1$, $H_2$, $H_3$ in general position in
$\tR P^8$. Fix a point $q\in \tR P^8-(H_1\cup H_2\cup H_3)$. Note that 
$X=\tR P^8-\{p\}$ is a non-orientable line bundle over $H_i$,
$i=1,2,3$. Choose sections $v_i$ of $X$ over $H_i$  such
that $\{v_i=0\}\cong\tR P^6$ meets $H_1\cap H_2\cap H_3$ transversally 
in $H_i$ for $i=1,2,3$. Now, $H_1\cap H_2\cap H_3\cong\tR P^5$. 
Let
$g\colon \tR P^5\to\tR P^8$ be the corresponding embedding.
Restricting $v_i$ to $\tR P^5$ we get three normal vector fields 
$v_1,v_2,v_3$ along $\tR P^5\subset \tR P^8$ and
$\{v_1=v_2=v_3=0\}\cong\tR P^2$.

Let $p\colon S^5\to\tR P^5$ be the universal cover. 
Then $h=g\circ p\colon S^5\to\tR P^8$ is an immersion. Let
$K_i=h^{-1}(\{v_i=0\})\cong S^4$. 
Choose Morse functions $\phi_i\colon S^5\to\tR$ such that
$\{\phi_i=0\}=K_i$. Let $\epsilon>0$ be 
small. Then $f\colon S^5\to\tR P^8$, given by
$$
f(x)=h(x)+\epsilon\sum_{i=1}^3\phi_i(x)v_i(h(x))\quad\text{for $x\in S^5$},
$$
is a generic immersion with $\Gamma_2(f)\cong\tR P^2$. Hence, $j_2(f)=1$.

\end{document}